\documentclass[12pt]{article}
\usepackage{amsmath,amsthm,amsfonts}
\usepackage{url}

\newcommand{\pf}{\noindent {\bf Proof:} }
\newtheorem{theorem}{Theorem}
\newtheorem{remark}{Remark}

\newtheorem{lemma}[theorem]{Lemma}

\title{Dejean's conjecture holds for $n\ge 27$}
\author{James Currie\thanks{The author is
supported by an NSERC Discovery Grant.} and
Narad Rampersad\thanks{The author is supported by an NSERC Postdoctoral
Fellowship.} \\
Department of Mathematics and Statistics \\
University of Winnipeg \\
515 Portage Avenue \\
Winnipeg, Manitoba R3B 2E9 (Canada) \\
\url{j.currie@uwinnipeg.ca} \\
\url{n.rampersad@uwinnipeg.ca}}

\begin{document}
\date{\today}
\maketitle

\begin{abstract}
We show that Dejean's conjecture
holds for $n\ge 27.$ This brings the final resolution of the conjecture by the approach of Moulin Ollagnier within range of the computationally feasible.
\end{abstract}

Repetitions in words have been studied since the beginning of the
previous century \cite{thueI,thue}. Recently, there has been much
interest in repetitions with fractional exponent
\cite{brandenburg,carpi,dejean,longfrac,krieger,mignosi}. For
rational $1<r\le 2$, a {\bf fractional $r$-power} is a non-empty
word $w=xx'$ such that $x'$ is the prefix of $x$ of length
$(r-1)|x|$. For example, $010$ is a $3/2$-power. A basic problem
is that of identifying the repetitive threshold for each alphabet
size $n>1$: \begin{quote}What is the infimum of $r$ such that an
infinite sequence on $n$ letters exists, not containing any factor
of exponent greater than $r$?\end{quote} The infimum is called the
{\bf repetitive threshold} of an $n$-letter alphabet, denoted by
$RT(n)$. Dejean's conjecture \cite{dejean} is that
$$RT(n)=\left\{\begin{array}{ll}7/4,&n=3\\
7/5,&n=4\\
n/(n-1)&n\ne 3,4\end{array}\right.$$

Thue, Dejean and Pansiot, respectively \cite{thue,dejean,pansiot}
established the values $RT(2)$, $RT(3)$, $RT(4)$. Moulin Ollagnier
\cite{ollagnier} verified Dejean's conjecture for $5\le n\le 11$,
and Mohammad-Noori and Currie \cite{morteza} proved the conjecture
for $12\le n\le 14$.

Recently, Carpi \cite{carpi} showed that Dejean's conjecture holds for
$n\ge 33$. Carpi's result is computation-free, and resolving Dejean's conjecture is thus reduced to filling a finite gap.
Conceptually, one would hope that the gap could now be filled from below, using the methods of \cite{ollagnier,morteza}. Since these approaches are
computationally intensive, optimizing Carpi's result is important. The present authors improved part of Carpi's
constructions to show that Dejean's conjecture holds for $n\ge
30$. (See \cite{archiv}.) In the present note we show that in fact
Dejean's conjecture holds for $n\ge 27$.

\begin{remark} Some months after the first draft of this paper, its goal has been vindicated: The final resolution of the conjecture via methods of Moulin Ollagnier becomes computationally feasible; in a recent paper the present authors proved Dejean's conjecture
by resolving computationally the cases $n\le 26$. Dejean's conjecture is correct! (See \cite{archiv2}.)
\end{remark}

The following definitions are from \cite{carpi}: For any
non-negative integer $r$ let $A_r=\{1,2,\ldots, r\}$. Fix $n\ge
27$. Let $m=\lfloor(n-3)/6\rfloor$. Let ker $\psi=\{v\in A_m^*
\mid \forall a\in A_m$, 4 {\it divides }$|v|_a\}.$ (We use this as
a definition; it is in fact the assertion of Carpi's Lemma~9.1.) A
word $v\in A_m^+$ is a $\psi$-{\bf kernel repetition} if it has
period $q$ and a prefix $v'$ of length $q$ such that $v'\in$ ker
$\psi$, $(n-1)(|v|+1)\ge nq-3.$ In \cite{archiv} we introduced the
following definition: If $v$ has period $q$ and its prefix $v'$ of
length $q$ is in ker $\psi$, we say that $q$ is a {\bf kernel
period} of $v$.

Let $B=\{0,1\}$ and let $S_n$ be the permutation group on $n$
elements. Consider the morphism $\phi:B^*\rightarrow S_n$
generated by
\begin{eqnarray*}\phi(0)&=&\begin{array}{lllllll}(1&2&3&&{\cdots} &(n-1))\end{array}\\
\phi(1)&=&\begin{array}{lllllll}(1&2&3&\cdots &(n-1)&n)\end{array}\end{eqnarray*}

This map is due to Pansiot \cite{pansiot}. A word $u\in B^*$ is a
{\bf $k$-stabilizing word} if $\phi(u)$ fixes $\{1, 2, 3, \ldots,
k\}$. The set of $k$-stabilizing words (for fixed $n$) is denoted by {\bf
$\mbox{Stab}_n(k)$.} Note that if $i<j$ then $\mbox{Stab}_n(j)\subseteq\mbox{Stab}_n(i).$

A map $\gamma_n:B^*\rightarrow A_n^*$ is defined by
$$\gamma_n(b_1b_2\cdots b_\ell)=a_1a_2\cdots a_\ell$$
where $a_i\phi(b_1b_2\cdots b_\ell)=1$ for $1\le i\le \ell$.

Carpi introduces a morphism $f:A_m^*\rightarrow B^*$ generated by

\begin{eqnarray*}
f(1)&=&y^px(101)^{2m}\\
f(a)&=&y^px(101)^{2m-2a}010(101)^{2a-1}
\end{eqnarray*}

where $2 \leq a \leq m$, $p=\lfloor n/2\rfloor$, $y$ is the suffix of $(01)^n$ of
length $n-1$ and $x$ is the suffix of $y$ of length $|y|-6m$.

The concepts of so-called {\bf short repetitions} and {\bf kernel
repetitions} were introduced by Moulin Ollagnier \cite{ollagnier}.
His work is complicated by the fact that his short repetitions are
words over $A_n$, while his kernel repetitions  are words over $B$ (although they code words
over $A_n$ via Pansiot's map).
Without going into the details, we recall that he reduced the
construction of an infinite word over $n$ letters attaining
threshold $n/(n-1)$ to avoiding
both short repetitions and kernel repetitions. Moulin Ollagnier's binary words
were fixed points of morphisms. In \cite{morteza}, a technique was
introduced for dealing separately with short repetitions and
kernel repetitions; the binary words given there can be viewed as
being produced by HD0L's: they have the form $g(h^\omega(0))$
where all words coded by $g(B^*)$ avoid short repetitions, and each $h$
is chosen to eliminate kernel repetitions.

Carpi's work follows essentially this strategy. The lemmas of his
paper show that $f(B^*)$ avoids short repetitions if $n\ge 30$.
For $m=5$ (corresponding to $n\ge 33$) he produces an infinite
word $w_5$ over $A_m$ such that $f(w_5)$ avoids kernel
repetitions. The exact statement of this division of work into short vs. kernel repetitions
 is in
his Proposition~3.2:
\begin{quotation}
{\bf Proposition 3.2:} Let $v\in B^*$. If a factor of
$\gamma_n(v)$ has exponent larger than $n/(n-1)$, then $v$ has a
factor $u$ satisfying one of the following conditions:
\begin{enumerate}
\item[(i)] $u\in \mbox{Stab}_n(k$) and $0 < |u| < k(n-1)$ for some $k\le n-1$
\item[(ii)] $u$ is a kernel repetition of order $n$.
\end{enumerate}
\end{quotation}

In our previous note, we improved only the second part of Carpi's
construction; he had shown that for $n\ge 30$, no factor $u$ of
$f(A_m^*)$ satisfied condition (i) above. As Carpi therefore
states at the beginning of section 9 of \cite{carpi}:
\begin{quotation}
By the results of the previous sections, at least in the case
$n\ge 30$, in order to construct an infinite word on $n$ letters
avoiding factors of any exponent larger than $n/(n-1)$, it is
sufficient to find an infinite word $w$ on the alphabet $A_m$ avoiding
$\psi$-kernel repetitions.
\end{quotation}
For $m=5$, Carpi was able to produce such an infinite word, based
on a paper-folding construction. He thus established Dejean's
conjecture for $n\ge 33$. The present authors refined this by constructing an infinite word $w_4$ on the
alphabet $A_4$ avoiding $\psi$-kernel repetitions. This
established Dejean's conjecture for $n\ge 30$. We remark that for $30\le n\le 32$ the word on
$A_n$ verifying Dejean's conjecture for $n$ is $\gamma_n(v)$, where
$v=f(w_4)$.

In the present note, we improve on the first aspect of Carpi's attack,
by showing that for $27 \le n\le 29$, no factor $u$ of $v=f(w_4)$
satisfies (i) above. This implies that Dejean's conjecture holds
for $n\ge 27$. Since $f$ is $r$-uniform where $r=(p+1)(n-1)$, to show that (i)
holds for $v$ it suffices to check that no factor $u\in f(B^3)$ satisfies (i).
In principle, this involves considering all factors of $f(B^3)$ of length less
than $(n-1)^2$. However, we shorten this computation considerably by combining several of
Carpi's lemmas.

\begin{lemma}
Suppose $n\ge 18$. Suppose that $u\in f(A_m^*)\cap \mbox{Stab}_n(k)$ and $|u|< k(n-1)$, some $k\in\{1,2,\ldots, n-1\}$.
Then $|u|=r(n-1)$ for some $r$, $p+1\le r<k\le 16$.
\end{lemma}
\pf Propositions and lemmas referenced in this proof are in \cite{carpi}. By Proposition~5.1, $k\ge 4$ so that $u\in\mbox{Stab}_n(4)$. It then follows from Proposition~6.3 that
$|u|\ge(p+1)(n-1)$. Since $|u|<k(n-1)$, we deduce that $k>p+1$. From $n\ge 18$ this means that $k>10$, so that surely $u\in\mbox{Stab}_n(7)$. Applying Lemma~7.1, we see that $|u|$ is divisible by $n-1$. We may thus write $|u|=r(n-1)$, $p+1\le r<k$. By the contrapositive of Proposition~7.2, $u\not\in\mbox{Stab}_n(17)$. It follows that $k\le 16.\Box$

We verify that Dejean's conjecture holds for $n=27,28,29$ by exhaustively examining factors $u$ of $f(B^3)$ of length $r(n-1)$ for $p+1\le r\le 15$, and verifying that such $u$ are not in $\mbox{Stab}_n(k)$ for any $k$, $r<k\le 16$. For $n=28,29$, the check only involves $r=15$, $k=16$. For $n=27$, we also must consider $r=14$. Code written in SAGE running on a PC performed the necessary verifications in about half an hour. The code is available at \begin{verbatim}www.uwinnipeg.ca/~currie/kstab.sage\end{verbatim}

\end{document}